\documentstyle[12pt]{article}    
\newcounter{props}[section]

\newcommand{\mda}{\rule[2.7pt]{1.0mm}{0.1mm}\hspace{1pt} }

\newcommand{\wrt}{with respect to\ }

\newcommand{\tx}{there exist\ }
\newcommand{\txs}{there exists\ }
\newcommand{\st}{such that\ }

\newcommand{\fe}{ for every\  }
\newcommand{\fa}{for each\  }

\newcommand{\sq}{sequence }

\newcommand{\isq}{increasing sequence\  }

\newcommand{\lra}{\longrightarrow }

\newcommand{\ra}{\rightarrow }

\newcommand{\calm}{\mbox{$\cal{M}$}}
\newcommand{\caln}{\mbox{$\cal{N}$}}

\newcommand{\cals}{\mbox{$\cal{S}$}}

\newcommand{\nin}{\mbox{$n\in\integers$}}

\newcommand{\calf}{\mbox{$\cal{F}$}}

\newcommand{\integers}{\mbox{${\bf N}$}}
\newcommand{\NN} {\mbox{${\bf N}$}}

\newcommand{\sumin}{\sum_{i=1}^{n}}

\newcounter{excounter}
\begin{document}           
\newcommand{\laga} {\left\{\rule{0mm}{4mm}\right.}
\newcommand{\lagb} {\left\{\rule{0mm}{5mm}\right.}
\newcommand{\lagc} {\left\{\rule{0mm}{6mm}\right.}
\newcommand{\lagd} {\left\{\rule{0mm}{7mm}\right.}
\newcommand{\raga} {\left.\rule{0mm}{4mm}\right\}}
\newcommand{\ragb} {\left.\rule{0mm}{5mm}\right\}}
\newcommand{\ragc} {\left.\rule{0mm}{6mm}\right\}}
\newcommand{\ragd} {\left.\rule{0mm}{7mm}\right\}}
\newcommand{\spp}{ {\rm supp}\,}
\newcommand{\sprt}{\vspace{0.2in}}
\newcommand{\sprth}{\vspace{0.1in}}
\newcommand{\sprtq}{\vspace{0.05in}}
\newcommand{\arrsp}{\vspace{2.5mm} \\}
\newcommand{\beq}{\begin{eqnarray*}}
\newcommand{\eeq}{\end{eqnarray*}}
\newcommand{\iffp}{it follows from Proposition\ }
\newcommand{\Iffp}{It follows from Proposition\ }
\newcommand{\hardexercise}{
\hspace{-2.5mm}$ ^{*}$\hspace{1.5mm}}
\newcommand{\sbs}{\subset}
\newcommand{\sps}{\supset}
\newcommand{\stm}{\setminus}
\newcommand{\mbf}{\mbox{ }\hfill}
\newcommand{\fmb}{\hfill\mbox{ }}
\newcommand{\proof}{
  \sprt \bf Proof.  \rm \  }
\newcommand{\prooft}{
  \sprt \bf Proof of Theorems 9.17 and 9.18. \noindent \rm \  }
\newcommand{\thm}[1]{\stepcounter{props}
  \sprt\bf \noindent \thesection .\theprops .\ 
  #1:\rm \  }
\newcommand{\sol}[1]{
  \sprt\hspace{1cm}\bf  #1\par\hspace{-5pt}\nopagebreak\rm \  }
\newcommand{\exercises}{\sprt
\[ {\rm EXERCISES}\]\sprt\small
\setcounter{excounter}{0}\begin{list}%
{\thesection --\theexcounter}{\usecounter{excounter}
}}
\newcommand{\abst}{\sprt
\[ {\rm ABSTRACT}\]
\begin{description}
\rightmargin=20 mm
}
\newcommand{\stareq}[3]{\vspace{7pt}$(#1)\hspace{#2cm}{\displaystyle
#3}$ \vspace{7pt}\\}
\newcommand{\bb}{
                 \item [ }
\Large
\[ {\bf Examples\; of\; Asymptotic}\; \ell_{1}\]

\[{\bf  Banach\; spaces}\]

\normalsize
\[{\rm by}\]
\[ {\rm S.\; A.\; Argyros\;\; and \;\; I.\; Deliyanni}\]
\[ {\rm Athens - Herakleion,\; Greece}\]

\abst
\item [\ \ \ \ \ \ \ \ ] {\em Two examples of 
asymptotic $\ell_{1}$ Banach spaces are given.
The first, $X_{u}$, has an unconditional basis and 
 is arbitrarily distortable. The second, $X$, does not
contain any unconditional basic sequence. Both are
spaces of the type of Tsirelson. We thus answer 
a question raised  by  W.T.Gowers.}
\end{description}


\mbf {\bf Introduction}\fmb

\sprth
The first example of an arbitrarily distortable Banach space
was constructed by Th. Schlumprecht in [Schl]. 
Schlumprecht's space was the starting point for the construction
by W.T. Gowers and B. Maurey of a Banach space not containing
an unconditional basic \sq (u.b.s.) [G-M]
and for the examples,  due  to  W.T. Gowers,
of a Banach space not containing $\ell_1$, $c_0$ or a
 reflexive subspace [G1] and of a space without u.b.s. but with
an asymptotically unconditional basis [G2].

A rapid development of the theory of Banach spaces followed 
the examples of Schlumprecht and Gowers -  Maurey.
We mention some results.

The notion of a hereditarily indecomposable Banach space 
was introduced in [G-M] and a new dichotomy property
for Banach spaces regarding this notion was proved by Gowers
[G3].  The remarkable solution of the distortion problem
for $\ell_p$, $(1<p< \infty )$ by E. Odell and
Th. Schlumprecht also makes use of  Schlumprecht's space.
Finally these results led to a new interest in the asymptotic
structure of Banach spaces [Mi-To], [Ma-Mi-To].

The examples we give in the present paper have as starting point
Tsirel\-son's celebrated example of the reflexive Banach space $T$
not containing any $\ell_p$.
We recall, following T. Figiel and W. Johnson [F-J] the definition
of Tsirelson's norm. Let $0<\theta <1$.
On $c_{00}$ (the space of finitely supported sequences)
we define implicitly the norm $\|\cdot\|_T$ by
\[ \| x\|_T =\max\left\{\| x\|_\infty ,\,
\sup\theta\sumin\| E_i x\|_T\right\} ,\]
where the $``\sup ''$ is taken over all families
$\{ E_1,E_2,\ldots ,E_k\}$ of finite subsets
of $\NN$ \st $n\leq E_1< E_2< \cdots < E_n$.
Tsirelson's space is an asymptotic $\ell_1$ space.
We recall the definition of this notion,
introduced in [Mi-To].

A Banach space with a normalized basis $\{ e_{k}\}_{k=1}^\infty$
is asymptotic $\ell_p$ if \txs a constant $C$ \st
\fe $n$ \txs $N=N(n)$ \st every \sq $(x_i)_{i=1}^n$ of
successive normalized blocks of $\{ e_{k}\}_{k=1}^\infty$ with
$N<\spp x_1<\spp x_2<\cdots <\spp x_n$ is $C$-equivalent
to the canonical basis of $\ell^n_p$.

We consider the following generalization of 
Tsirelson's example.
Let $\calm$ be a family of finite subsets of $\NN$ closed
in the topology of pointwise convergence.
A finite \sq $\{ E_i\}_{i=1}^n$ of finite subsets
of $\NN$ is said to be $\calm$-admissible if \txs
 a set $F=\{ k_1,\ldots ,k_n\}\in\calm$ \st
\[ k_1\leq E_1<k_2\leq E_2<\cdots k_n\leq E_n .\]
Let $0<\theta <1$.  The Tsirelson type Banach space
$T[\calm ,\theta ]$ is  the completion of $c_{00}$ under the
norm $\| \cdot\|_{{\cal M} ,\theta}$ which is defined by the 
following implicit  equation:
\[ \| x\|_{{\cal M} ,\theta} =\max\left\{\| x\|_\infty ,\,
\sup\theta\sumin\| E_i x\|_{{\cal M} ,\theta}\right\} ,\]
where the $``\sup ''$ is taken over all   $n$ and all
$\calm$-admissible sequences $\{ E_i\}_{i=1}^n$.
It is clear that Tsirelson's original space is 
$T[\cals ,\theta ]$ where $\cals$ is the Schreier family
defined by
\[ \cals =\{ F: F\sbs\NN, \# F\leq \min F\} .\]
Consider $A_n=\{ F: F\sbs\NN, \# F\leq n\} .$
S. Bellenot, [B], has proved the following result:
For every $1<p<\infty$ and $n\geq 2$ \txs $0<\theta <1$ 
\st $T[A_n ,\theta ]$ is isomorphic to $\ell_{p}$.
The spaces $T[\calf_\xi ,\theta ]$ (
the families  $\calf _\xi$, $\xi <\omega_1$ are defined below) were
introduced by the first 
author in order to prove  the following result:
For every  $\xi <\omega_1$ \txs a reflexive Banach space
$T_\xi$ \st every infinite dimensional subspace
of $T_\xi$ has Szlenk index greater than $\xi$
(preprint, 1987). 
The general spaces $T[\calm ,\theta ]$ were defined in
[Ar-D].

The examples we present here are defined using 
``mixed Tsi\-rel\-son's norms''. These norms are defined 
by sequences $\{ \calm_n\}_{n=1}^\infty$ and 
$\{ \theta_n\}_{n=1}^\infty$ \st each $\calm_n$ is a
family  of finite subsets of $\NN$ closed in the 
topology of pointwise convergence and $0<\theta_n <1$,
$\lim_{n\ra\infty}\theta_n =0$.
The norm in the space $T[(\calm_n,\theta_n)_{n=1}^\infty ]$ is
defined by
\[ \| x\| =\max\left\{\| x\|_\infty ,\,
\sup_k\left\{ \theta_k \sup\sumin\| E_i x\|\right\}\right\} ,\]
where the inner  $`\sup '$ is taken over all   $n$ and all
$\calm_k$-admissible families $( E_1,\ldots ,E_n)$.
It is easy to see that if the Schreier  family $\cals$ is
contained in one of the families $\calm_n$ then the space
$T[(\calm_n ,\theta_n )_{n=0}^\infty]$ is
asymptotic-$\ell^1$. 
Our first space, $X_u$,  is a space of the form 
$T[(\calm_n ,\theta_n )_{n=0}^\infty]$
for appropriate sequences $(\calm_n)_{n=1}^\infty$ and
$(\theta _n)_{n=1}^\infty$. $X_u$ has an unconditional basis
and is arbitrarily distortable. 
The second  space, $X$, does not contain any unconditional basic 
sequence. In fact it is hereditarily indecomposable. 
$X$ is constructed via $X_u$ in a way
similar to the one used in [G-M] to pass from 
Schlumprecht's space to the Gowers-Maurey space.
The basic idea for this comes from the fundamental construction
by Maurey and Rosenthal [M-R] of a weakly null \sq
without unconditional basic subsequence.

Although our approach is different from that of 
Sclumprecht, Gowers and Maurey, it seems that 
the ingredients needed  for the proofs are similar.
So, for example, the normalized $(\epsilon ,j)$-special
convex combinations correspond to $\ell^1_N$ vectors
and the rapidly increasing $(\epsilon ,j)$-s.c.c.'s
correspond to sums of rapidly increasing sequences.

\newpage
\section{ Preliminaries} 

\sprt
{\bf (a) Generalized Schreier families}

\sprth
The {\em Schreier family} \cals\ is the set of all finite
subsets of \NN\ satisfying the property $\# A\leq \min A$.
It is easy to see that this family is closed in the topology
of pointwise convergence.

\sprt
{\bf 1.1 Definition.} Given $\calm ,
\caln$, families of finite subsets of \NN\ 
which are closed in the topology of pointwise convergence,
the $\calm$ {\em operation on} $\caln$ is defined as
\beq \calm [\caln ]\!\!\!\!&=&\!\!\!
\left\{ F\sbs\NN : F=\bigcup_{i=1}^{s}F_{i},\right.
\;\; s\in \NN ,\;
 m_{1}\leq F_{1}<m_{2}\leq F_{2}<\cdots <m_{s}\leq F_{s}\\
\!\!\!& &\!\!\!\left.\;\;\;\;\;\;\;\;\;\;\;\;\;\;\;\;
\; F_{i}\in \caln , i=1,\ldots ,s\; {\rm and}\;\{ m_{1},\ldots , m_{s}\}
\in \calm\rule{0mm}{6.5mm}\right\} .\eeq
$\calm [\caln ]$ is a family of finite subsets of \NN\ which is closed
in the topology of pointwise convergence.

\sprt
{\bf 1.2 Definition.} The generalized Schreier families
$\{\calf_{\xi}\}_{\xi <\omega_{1}}$
are defined as follows:

 $\calf_{0}=\{\{ n\} :n\in\NN\}$,

$\calf_{\xi +1}=\cals [\calf_{\xi}]$. 

For $\xi$  a limit
ordinal we let $\{\xi_{n}\}_{n=1}^{\infty}$ be a fixed  \sq 
strictly increasing to $\xi$ and set
\[ \calf_{\xi}=\left\{ A\sbs\NN :n\leq\min A\;{\rm and}\;
A\in\calf_{\xi_{n}}\right\} .\]
The families $\{\calf_{\xi}\}_{\xi <\omega_{1}}$ have been 
introduced in [Al-Ar].

{\bf Remark.}
It is easy to see that for $\xi_1,\xi_2$ \tx 
$\xi<\omega_1$ \st $\calf_{\xi_1}[\calf_{\xi_2}]
\sbs \calf_{\xi}$. In particular, for $m,n\in\NN$,
$\calf_n[\calf_m] =\calf_{n+m}$.

\sprt
{\bf (b) Tsirelson type spaces}

\sprth
 In [Ar-D] a space $T[\calm ,\theta ]$ has been
defined, where \calm\ is a family of finite subsets of \NN\
closed in the topology of pointwise convergence and
$\theta$ a real number with $0<\theta <1$.

We recall that definition. Given $\calm$ as above, a family
$(E_{1},\ldots ,E_{n})$ of succesive finite subsets of \NN\
is said to be $\calm${\em -addmissible} if \txs a set
$A=\{ m_{1},\ldots , m_{n}\}\in\calm$ \st
$m_{1}\leq E_{1}<m_{2}\leq E_{2}<\cdots <m_{n}\leq E_{n}$.
The norm on the space $T[\calm,\theta]$ is defined 
implicitly by the formula 
$\| x\| =\max\{\|x\|_{\infty},\theta\sup\sum_{i=1}^{n}\|
E_{i}x\|\}$, where the `sup' is taken over all $n$ and
all \calm -admissible 
$(E_{1},\ldots ,E_{n})$.

It is known that  if the Cantor-Bendixson
index of \calm\ is greater than $\omega$, then the space
$T[\calm,\theta]$ is reflexive. In 1.3 we prove a somewhat more
general result.

\sprt
{\bf (c) Mixed Tsirelson norms}

\sprth
Let $\{ \calm_k\}_{k=1}^{\infty}$ be families of finite subsets
of $\NN$ \st for each $k$:

\newcounter{lista}
\begin{list}{(\alph{lista})} {\usecounter{lista}}
\item $\calm_{k}$ is closed in the topology of pointwise
convergence.
\item $\calm_k$ is adequate, i.e. if $A\in\calm$ and
$B\sbs A$ then $B\in\calm_k$.
\item The Cantor - Bendixson index of $\calm_k$ is greater
than $\omega$.
\end{list}
Let $\{\theta_k\}_{k=1}^{\infty}$ be a \sq of positive
reals with each $\theta_k<1$ and $\lim \theta_k =0$.\\
Then the mixed Tsirelson norm defined by 
$(\calm_k ,\theta_k)_{k=1}^{\infty}$ is given by the implicit
relation
\[ \| x\| =\max\left\{ \| x\|_\infty,\sup_{k}\left\{\theta_k
\sup\sumin\| E_i x\|\right\}\right\} ,\]
where the inside $``\sup ''$ is taken over all
$\calm_k$-admissible families $E_1,\ldots ,E_n$.

The Banach space defined by this norm is denoted by 
$T\left[ (\calm_k,\theta_k )_{k=1}^{\infty}\right]$.

\sprt
{\bf 1.3 Proposition.}
The space $X=T\left[ (\calm_k,\theta_k )_{k=1}^{\infty}\right]$
is reflexive and $\{ e_n\}_{n=1}^{\infty}$ is an 1-unconditional
basis for $X$.

\proof
The proof is similar to the original proof of Tsirelson in [T].
We first give an alternative definition of the norm of $X$.

We define inductively the following sets:

$K^0=\{ \lambda e_n:\nin , |\lambda |\leq 1\} .$\\
Given $K^s$,

$K^{s+1}=K^s\cup \laga\theta_k (f_1+\cdots +f_d\}):
k\in\NN , d\in\NN , f_i\in K^s , i=1,\ldots ,d,$\\
$\;\;\;\;\;\;\;\spp f_1<\spp f_2<\cdots <\spp f_d$ and the set
$\{ \spp f_1,\ldots ,\spp f_d\}$ is 
$\calm_k$-admissible $\raga .$

Finally, we set
\[ K=\bigcup_{s=0}^{\infty}K^s\]
and for $x\in c_{00}$ we define
\[ \| x\|=\sup_{f\in K}<x,f>.\]
Then $X$ is the completion of $(c_{00},\|\cdot\| )$.
It is easy to see that $\{ e_n\}_{n=1}^{\infty}$ is a
1-unconditional basis for $X$.

To show that $X$ is reflexive, we have to show that the basis
$\{ e_n\}_{n=1}^{\infty}$ is shrinking and boundedly complete.

(a) $\{ e_n\}_{n=1}^{\infty}$ is a shrinking basis for $X$.

Let $\theta=\max_k\theta_k<1$. For $f\in X^*$ and $m\in\NN$,
denote by $Q_m(f)$ the restriction of $f$ to the space
generated by $\{ e_k\}_{k\geq m}$. It suffices to prove the 
following: For every $f\in B_{X^*}$ there is $m\in\NN$
\st $Q_m(f)\in \theta B_{X^*}$. Recall that $B_{X^*}=
\overline{co(K)}$ where the closure is in the topology of pointwise 
convergence.
We shall first prove the following: 

{\em Claim.} For every $f\in \overline{K}$
there  is $m$ \st $Q_m(f)\in \theta\overline{K}$.

So, let $f\in \overline{K}$ and let $\{ f^n\}_{n=1}^\infty$
be a \sq in $K$ converging pointwise to $f$.

If $f^n\in K^0$ for an infinite number of $n$, we have 
nothing to prove.
So suppose that \fe $n$ there are $k_n\in\NN$, a set
$\{ m^n_1,\ldots ,m^n_{d_n}\}\in \calm_{k_n}$ and vectors
$f_i\in K$, $i=1,\ldots ,d_n$ \st
$m^n_1\leq \spp f^n_1<
m^n_2\leq \spp f^n_2<\cdots <
m^n_{d_n}\leq \spp f^n_{d_n}$.
If there is a sub\sq of $\theta _{k_n}$ converging to 0,
then $f=0$. So we may suppose that there is a $k$ \st
$k_n=k$ for all $n$, i.e. $\theta _{k_n}=\theta_{k}$ and
$\{ m^n_1,\ldots ,m^n_{k}\}\in \calm_{k}$.

Since $\calm_k$ is compact, substituting $f^n$ with a sub\sq
we get that there is  a set
$\{ m_1,\ldots ,m_d\}\in\calm_k$ \st the \sq of indicator
functions of the sets 
$\{ m^n_1,\ldots ,m^n_{d_n}\}$ converges to the indicator function\\
of $\{ m_1,\ldots ,m_d\}$.
So, for large $n$,
$m^n_i=m_i$, $i=1,\dots ,d$, and $m^n_{d+1}\ra \infty$ as
$n\ra\infty$.\\
Passing to a further sub\sq of $(f^n)_{n=1}^\infty$, we get
that there exist $f_i\in\overline{K}$,
$i=1,\ldots ,d$ with $\spp f_i\sbs [m_i,m_{i+1})$,
$i=1,\ldots ,d-1$ and $\spp f_d\sbs [m_d,\infty )$ \st
$f^n_j\ra f_j$ pointwise for
$j=1,\ldots ,d$. 
We conclude that $f=\theta_k(f_1+\cdots +f_d)$,
so $Q_{m_d}(f)=\theta_kf_d\in \theta\overline{K}$.

The proof of the claim is complete.
In particular we get that $\overline{K}$ is a weakly
 compact subset of $c_0$.

By standard arguments we can now pass to the case of
$B_{X^*}=\overline{co(K)}$.

\sprtq
(b) $\{ e_n\}_{n=1}^\infty$ is a boundedly complete
basis for $X$.

Suppose on the contrary that \tx $\epsilon >0$ and a
block \sq $\{ x_i\}_{i=1}^\infty$ of 
 $\{ e_n\}_{n=1}^\infty$ \st $\|\sum_{i=1}^\infty x_i\|\leq 1$
while $\| x_i\|\geq \epsilon $ for $i=1,2,\ldots$.

Choose $n_0\in \NN$ \st $n_0\theta_1>\frac{1}{\epsilon }$.
Using the fact that the $n_0+1$-derivative set of 
$\calm_1$ is non-empty, one can choose a set
$\{ m_1,\ldots ,m_{n_0}\}\in \calm_1$ and a subset
$\{ x_{i_k}\}_{k=1}^{n_0}$ of $\{ x_i\}_{i=1}^\infty$  \st
\[ m_1\leq\spp x_{i_1}<m_2\leq\spp x_{i_2}<\cdots <
m_{n_0}\leq\spp x_{i_{n_0}} .\]
Then
\[ \left\|\sum_{k=1}^{n_0}x_{i_k}\right\|\geq
\theta_1\sum_{k=1}^{n_0}\| x_{i_k}\geq n_0\theta_1>1,\]
a contradiction and the proof is complete.

\sprt
{\bf (d) Special convex combinations}

\sprth
Next we prove a property of $\{\calf_n\}_{n=1}^\infty$
which is important for our constructions.

\sprt
{\bf 1.4 Proposition.} For every $n\in \NN$, $\epsilon >0$,
\txs $m>n$ \st \fe  infinite subset $D$ of $\NN$ \txs 
a convex combination $x=\sum_{n=1}^{s}a_ne_{k_n}$
with $\spp x\sbs D$,
$\spp x\in \calf_m$, $\{ a_n\}_{n=1}^s$ is in decreasing
order and $|x|_n<\epsilon $, where $|\cdot |_n$ denotes the norm of
the space $T\left[\calf_n ,\frac{1}{2}\right]$.

\sprth
We prove first the following

\sprt
{\bf 1.5 Lemma.} For every $n\geq 2$, $\epsilon >0$, $D$ infinite
subset of $\NN$, \txs a convex combination 
$x=\sum_{n=1}^{s}a_ne_{k_n}$ \st
\begin{list}{(\roman{lista})} {\usecounter{lista}}
\item $\spp x\sbs D$, $\spp x\in \calf_n$
\item $\{ a_n\}_{n=1}^s$ is in decreasing order
\item For every $F$ in $\calf_{n-1}$ we have that $\sum_{n\in F}
a_n<\epsilon $.
\end{list}

\proof
We prove it for $n=2$. The general case is proved by induction
in a similar manner.

Choose $n_0$ \st $\frac{1}{n_0}<\epsilon ^2$, $n_0\in D$.
Set $A_1=\{ n_0\}$ and choose $A_2$ in
$\calf_1$ \st $A_2\sbs D$, $A_1<A_2$, $\# A_2\geq
\frac{n_0}{\epsilon ^2}$.
Similarly choose $A_3,\ldots ,A_{n_0}$. It is easy to check
that the convex combination
\[ \frac{1}{n_0}\sum_{l=1}^{n_0}\frac{1}{\# A_l}\cdot
\chi_{A_l}\]
is the desired vector. The proof is complete.

\sprth
{\bf Proof of Proposition.}
Choose $l$ \st $\frac{1}{2^l}<\frac{\epsilon }{2}$ and $m$ \st
$\calf_n[\cdots [\calf_n]\cdots ]$ ($l$-times) is
contained in $\calf_{m-1}$.

Next  choose, by the previous lemma, a convex combination
satisfying (i), (ii) and (iii) for the given $D$ and
$\frac{\epsilon }{2}$. We claim that this is the desired vector.

Indeed, choose any $\phi$ in the norming set $K$, as it is defined
in the proof of Proposition 1.3 for the space
$T \left[ \calf_n,\frac{1}{2}\right]$. Notice that
the set $L=\left\{ k\in\NN_r:|\phi (x)|\geq\frac{1}{2^l }\right\}$
is in $\calf_{m-1}$. To see this one proves that
$\phi |L$ belongs to $K^l$ and \fe $\phi '$ in $K^l$,
$\spp \phi '$ is in $\calf_n[\cdots [\calf_n]\cdots ]$ ($l$-times),
hence indeed $\spp (\phi|L)$ is in $\calf_{m-1}$. Therefore
\[ \left|\rule{0mm}{4mm}\phi (x)\right|\leq 
\left|(\phi|L)(x)\rule{0mm}{4mm}\right| +
\left|\rule{0mm}{4mm}(\phi |L^c)(x)\right|\leq\frac{\epsilon }{2}+
\frac{\epsilon }{2}=\epsilon .\]
The proof is complete.

\section{ Definition of the space $X_{u}$}

\sprth
Choose $\{ m_{j}\}_{j=0}^{\infty}$
\st $m_0=2$, and $m_j>m^{m_{j-1}}_{j-1}$.
Inductively we choose a family $\{ \calm_j\}_{j=0}^\infty$
\st each $\calm_j$ is a family of finite subsets of $\NN$ closed
in the topology of pointwise convergence.

We set $\calm_0=\calf_1$.

Suppose that $\{ \calm_j\}_{j=1}^n$ has been constructed
so that it is increasing and \fe $j$ \txs $k_j$ \st
$\calm_j\sbs \calf_{k_j}\sbs\calm_{j+1}$.
Choose $s_{n+1}$ \st for $\epsilon _{n+1}=\frac{1}{m^4_{n+2}}$,
$\calf_{s_{n+1}}$ satisfies the assumptions of Proposition 1.4
for $\calf_{k_n}$, $\epsilon _{n+1}$.
Choose $l_{n+1}$ \st $\frac{1}{2^{l_{n+1}}}<
\frac{1}{2m_{n+1}}$ and set
\[ \calm_{n+1}=\calf_{s_{n+1}}[\cdots [\calf_{s_{n+1}}]
\cdots ],\;\; l_{n+1}\;{\rm times}.\]
The space
\[ X_{u}=T\left[\left(\calm_{j},\frac{1}{m_{j}}
\right)_{j=0}^{\infty}\right]\]
is a reflexive Banach space with an unconditional basis.

For later use we need an explicit definition of the family of functionals
that  define the norm on the space $X_{u}$.

We set $K^{0}_{j}=\{\pm e_{n}:n\in\NN\}$.

Assume that $\{ K^{n}_{i}\}_{i=0}^{\infty}$ have been defined.
Then we set
\beq 
K^{n}&=&\bigcup_{i=0}^{\infty}K^{n}_{i}\\
K^{n+1}_{j}&=&K^{n}_{j}\cup\left\{\frac{1}{m_{j}}
(f_{1}+\cdots +f_{d}):\{ {\spp}f_{1}<\cdots <{\spp}f_{d}\}
\;{\rm is}\rule{0mm}{7mm}\right.\\
& &\;\;\;\;\;\;\;\;\;\;\;\;\;\;\;\;\;\;\;\;\;\;\;
\left. \calm_{j}{\rm \mda admissible\; and\;} f_{1},
\ldots ,f_{d}\;{\rm belong\; to}\;K^{n}\rule{0mm}{5mm}\right\} .\eeq
Set  $K=\cup_{n=0}^{\infty}K^{n}$.
Then the norm $\|\cdot\|$  on $X_{u}$ is 
\[ \| x\| =\sup\left\{ f(x):f\in K
\right\} .\]

\sprth{\bf Notation} For $j=0,1,\ldots$ we denote by
$K_{j}$ the set 
\[ K_{j}=\bigcup_{n=1}^{\infty}K^{n}_{j}.\]

We will need the following property of the families 
$\calm_{n}$ which can be easily proved by induction:

\sprt{\bf 2.1 Lemma.}
Let $n\in\NN$ and $F=\{ s_{1}<s_{2}<
\cdots <s_{d}\}\sbs\NN$ with $F\in\calm_{n}$.
If $G=\{ t_{1}<t_{2}<\cdots <t_{r}\}$ is \st
$r\leq d$ and
$s_{p}\leq t_{p}$ for $p=1,2,\ldots ,r$, then
$G\in \calm_{n}$.

\sprt
$\!${\bf 2.2 Notation.} We denote by $\| \cdot\|_{j}$
the norm of the space $T\left[\left(\calm_{n},
\frac{1}{m_{n}}\right)_{n=0}^{j}\right]$ and by $\|\cdot\|_{j}^{*}$
the corresponding dual norm.
Notice that since $\calm_j$ is a subfamily of $\calf_{k_j}$,
$\| x\|_{j}\leq |x|_{k_j}$ where $|\cdot |_{k_j}$ is the norm
of $T\left[\calf_{k_j},\frac{1}{2}\right] .$

\sprt
{\bf 2.3 Definition.} Let $m\in\NN$, $\phi\in 
K^{m}\stm K^{m-1}$.
We call {\em analysis} of $\phi$ any \sq $\{ K^{s}(\phi)
\}_{s=0}^{m}$ of subsets of $K$ 
such that:

\begin{list}{\arabic{lista})} {\usecounter{lista}}
\item For every $s$, $K^{s}(\phi )$ consists of successive
elements of $K^{s}$ and\linebreak
$\cup_{f\in K^{s}(\phi)}{\spp}f={\spp}\phi$.
\item If $f$ belongs to $K^{s+1}(\phi)$ then either
$f\in K^{s}(\phi )$ or \txs $j$ and
$f_{1},\ldots ,f_{d}\in K^{s}(\phi )$ with
$\{\spp f_{1}<\cdots <\spp f_{d}\}$ $\calm_{j}$-admissible
and \st $f=\frac{1}{m_{j}}(f_{1}+\cdots +f_{d})$.
\item $K^{m}(\phi )=\{\phi\}$.
\end{list}

\sprt
{\bf Remark.} Every $\phi\in K$
 has an analysis. Also, if
$f_{1}\in F^{s}(\phi)$, $f_{2}\in F^{s+1}(\phi )$, then
either ${\spp}f_{1}\sbs{\spp}f_{2}$ or
${\spp}f_{1}\cap{\spp}f_{2}=\emptyset$

\sprt
{\bf 2.4 Definition.} (a) Given $\phi\in K^{m}\stm
K^{m-1}$ and $\{ K^{s}(\phi )\}_{s=0}^{m}$ a fixed
analysis of $\phi$, then for a given finite block \sq 
$\{ x_{k}\}^{l}_{k=1}$ we set
\[ s_{k}=\left\{\begin{array}{l}
\max\{ s:0\leq s<m{\rm \;and\; there\; are\; at\; least\; two}\;
f_{1}, f_{2}\in F^{s}(\phi )\;{\rm such}\\
\;\;\;\;\;\;\;\;\;\;\;\;\;\; {\rm that}\;
{\spp}f_{i}\cap{\spp}x_{k}\neq\emptyset ,i=1,2\}\arrsp
0\;\;\;{\rm if}\; \#\spp x_{k}\leq 1\end{array}\right.\]

(b) For $k=1,\ldots ,l$, we define the {\em initial} and 
{\em final part\/}
of $x_{k}$ \wrt $\{K^{s}(\phi )\}_{s=0}^{m}$,
denoted by $x_{k}'$ and $x_{k}''$ respectively,
as follows: Let $\{ f\in F^{s_{k}}(\phi ):
\spp f\cap \spp x_{k}\neq\emptyset\} =
\{f_{1},\ldots ,f_{d}\}$, where $\spp f_{1}<
\cdots <\spp f_{d}\}$. Then we set 
$x_{k}'=x_{k}|\spp f_{1}$, 
$x_{k}''=x_{k}|\cup_{i=2}^{d}\spp f_{i}$.

\sprt
{\bf A. Estimates on the basis} $(e_{n})_{n\in{\bf N}}$

\sprth
{\bf 2.5 Definition.} Given $\epsilon >0$, an 
$(\epsilon,j)${\em -basic special convex combination} 
($(\epsilon ,j)$-s.c.c.) is an element of the form
$\sum_{k\in F}a_{k}e_{k}$ \st $F\in\calm_{j}$,
$ a_{k}\geq 0$, $\sum_{k\in F} a_{k}=1$ and
$\|\sum_{k\in F} a_{k}e_{k}\|_{j-1}<\epsilon$.

\sprth
{\bf Remark.}
Proposition 1.4 and the definition of $\calm_j$ guarantee 
the existence of $(\epsilon ,j)$-basic s.c.c.'s for
$\epsilon =\frac{1}{m^4_{j+1}}$ and with the additional property
that the coefficient $(a_k)_{k\in F}$ are in decreasing order.

So, when referring to an $(\epsilon ,j)$-basic s.c.c.
$\sum _{k\in F}a_{k}e_{k}$ we will always mean that the
$a_{k}$'s are in decreasing order.

\sprt
{\bf 2.6 Proposition.} For  given $j\in\NN$, $\epsilon <
\frac{1}{m^{3}_{j}}$ and $\sum_{k\in F}a_{k}e_{k}$ 
an $(\epsilon ,j)$-s.c.c. we have that: For 
$\phi\in K$
\[ \left|\phi\left(\sum_{k\in F}a_{k}e_{k}\right)\right|\leq
\left\{\begin{array}{ll}
\frac{1}{m_{s}} & {\rm if}\; \phi\in K_{s}, s\geq j\arrsp
\frac{2}{m_{s}\cdot m_{j}} & {\rm if}\; \phi\in K_{s}, s< j.
\end{array}\right.\]

\sprth
{\bf Proof:} If $s\geq j$ then the estimate is obvious.

Assume that $s<j$ and for some $\phi\in K_{s}$,
$|\phi (\sum  a_{k}e_{k})|>
\frac{2}{m_{s}m_{j}}$. 
Without loss of generality we assume that
$\phi (e_{k})\geq 0$ for all $k$.
Then 
$\phi =\frac{1}{m_{s}}(x^{*}_{1}+\cdots +x^{*}_{d})$,
where $\{ \spp x^{*}_{1}<\cdots <\spp x^{*}_{d}\}$ is
$\calm_{s}$-admissible.
We set
\[ D=\left\{ k\in F:\sum_{i=1}^{d}x^{*}_{i}(e_{k}) >
\frac{1}{m_{j}}\right\}.\]
Then $\sum_{k\in D}a_{k}>\frac{1}{m_{j}}$.
If not, then 

\noindent$ (\sum_{i=1}^{d}x^{*}_{i})
(\sum_{k\in F}a_{k}e_{k})\leq
 (\sum_{i=1}^{d}x^{*}_{i})
(\sum_{k\in D}a_{k}e_{k}) +
 (\sum_{i=1}^{d}x^{*}_{i})
(\sum_{k\not\in D}a_{k}e_{k}) \leq
\frac{2}{m_{j}}$, hence
$\phi (\sum_{k\in F} a_{k}e_{k})\leq
\frac{2}{m_{s}m_{j}}$, a contradiction.
Set $y^{*}_{i}=x^{*}_{i}|D$. Then it is easy to see that
$\left\|\frac{1}{m_{s}}(y^{*}_{1}+\cdots +y^{*}_{d})
\right\|_{j-1}^{*}\leq 1$ and 
$\frac{1}{m_{s}}(y^{*}_{1}+\cdots +y^{*}_{d})
\left(\sum_{k\in F} a_{k}e_{k}\right)\geq 
\frac{1}{m^{s}m_{j}^{2}}>\frac{1}{m^{3}_{j}}$,
 a contradiction and the proof
is complete.

\sprt
{\bf 2.7 Remark.} (a) It is easy to see that every $(\epsilon ,j)$-s.c.c.
 in $X_{u}$ has norm greater than  or equal to $\frac{1}{m_{j}}$.
Therefore, for $\epsilon <\frac{1}{m^{3}_{j}}$, we get that the
norm of the $(\epsilon ,j)$-basic s.c.c. is exactly $\frac{1}{m_{j}}$.

(b) It is crucial for the rest of the proof that for $s<j$,
and $x^{*}_{i}\in K$, $i=1,\ldots ,d$ with 
$\{ \spp x^{*}_{i}\}_{i=1}^{d}$ $\calm_{s}$-admiissible,
\[\left|\frac{1}{m_{s}}(x^{*}_{1}+\cdots +x^{*}_{d})\right|
\left(\sum_{k\in F} a_{k}e_{k}\right)\leq 
\frac{2}{m_{s}m_{j}}.\]
In other words, for the normalized vector 
$m_{j}\sum_{k\in F}a_{k}e_{k}$ we have that

$\displaystyle \left|\sum_{i=1}
^{d}x_{i}^{*}\right| \left(\left( m_{j}(\sum_{k\in F}a_{k}e_{k}
\right)\right)\leq 2$. 

\sprt
{\bf B. Estimates on block sequences}

\sprt
{\bf 2.8 Definition.} (a) Given a normalized block 
\sq $(x_{k})_{k\in{\bf N}}$ in $X_{u}$,
a convex combination $\sumin a_{i}x_{k_{i}}$ is said to
be an $(\epsilon ,j)${\em -s.c.c.} if 
\tx $l_1<l_2<\cdots l_n$ \st 
$\spp x_{k_1}\leq l_1<
\spp x_{k_2}\leq l_2<\cdots <
\spp x_{k_n}\leq l_n$ and
$\sum_{i=1}^{n} a_{i}e_{l_{i}}$ is an $(\epsilon ,j)$-basic s.c.c.

(b) An $(\epsilon ,j)$-s.c.c. is called {\em normalized\/}
if $\|\sumin a_{i}x_{k_{i}}\|\geq\frac{1}{2}$.

\sprt{\bf 2.9 Lemma.} Let $\{x_{k}\}_{k=1}^{\infty}$ be a
normalized block \sq in $X_{u}$ and 
$j=0,1,2,\ldots$, $\epsilon >0$, then \txs 
$\{ y_{k}\}_{k=1}^{n}$ a finite block \sq of
$\{x_{k}\}_{k=1}^{\infty}$ \st $\| y_{k}\| =1$
and a convex combination $\sum a_{k}y_{k}$ is an
$(\epsilon ,j)$-s.c.c. with
$\|\sum a_{k}y_{k}\|>\frac{1}{2}$.

\sprth {\bf Proof.}
Choose $\{y_{k}^{1}\}_{k=1}^{\infty}$ successive blocks of
$\{x_{k}\}_{k=1}^{\infty}$ \st each $y_{k}^{1}$ is 
an $(\epsilon ,j)$-s.c.c. of $\{x_{k}\}_{k=1}^{\infty}$
defined by an $(\epsilon ,j)$-basic s.c.c. $z^{1}_{k}$
\st
$\spp z^1_k\in \calf_{s_j}$, $k=1,2,\ldots$.
If for some $k_{0}$, $\| y^{1}_{k}\|\geq\frac{1}{2}$
then we are done; if not we consider the normalized block \sq
$x^{1}_{k}=\frac{y^1_k}{\| y^{1}_{k}\|}$ and 
apply the same procedure for 
$\{x_{k}^{1}\}_{k=1}^{\infty}$ as we did for 
$\{x_{k}\}_{k=1}^{\infty}$. 
So we get $\{y_{k}^{2}\}_{k=1}^{\infty}$, a \sq of
$(\epsilon ,j)$-s.c.c.'s \st
each $y^2_k$ is defined by a basic s.c.c. $z^2_k$ with
$\spp z^2_k\in\calf_{s_j}$. 
If $\| y^2_k\|\leq \frac{1}{2}$ then
$\{ \spp x_k : \spp x_k\sbs \spp y^2_{k_0}\}$ is
$\cals [\calf_{s_j}[\calf_{s_j}]]$-admissible
(so $\calm_j$-admissible), we get that
$\frac{1}{m_{j}}\leq\left\|\frac{1}{2} y^{2}_{k_{0}}\right\|
 <\frac{1}{2^{2}}$.

Repeating the procedure $l_{j}$ times, if  we never  get a
$y^{i}_{k}$, $1\leq i\leq l_{j}$, with  
$\| y_{k}^{i}\|\geq \frac{1}{2}$, then
we arrive at a $y^{l_{j}}_{k_{0}}$ \st
$\{\spp x_{k}:\spp x_{k}\sbs \spp y_{k_{0}}^{l_{j}}\}$
is $\cals [\calm_{j}]$-admissible and
\[ \frac{1}{2m_{j}}\leq\frac{1}{2^{l_{j}-1}}\left\|
y^{l_{j}}_{k_{0}}\right\| <\frac{1}{2^{l_{j}}},\]
a contradiction since $2m_{j}<2^{l_{j}}$ and the proof
is complete.

\sprth The way of proving the above result is similar to
the one of Gowers and Maurey [G-M] on the existence
of $\ell^{1}_{N}$-vectors.


\sprt
{\bf 2.10 Proposition.} Let $j\in\NN$. Let 
$\{ x_{k}\}_{k=1}^{n}$ be a finite block \sq of normalized
vectors in $X_{u}$.
Let $\{ l_1,\ldots ,l_n\}$ be \st
$\spp x_{k_1}\leq l_1<
\spp x_{k_2}\leq l_2<\cdots <
\spp x_{k_n}\leq l_n$ and suppose that
 $\{ l_1,\ldots ,l_n\}\in \calm_j$.
Then, for every
$q\leq j$ and every $\phi\in K_{q}$, \txs $\psi\in co(K_{q})$
\st $|\phi (x_{k})|\leq 2\psi (m_{j}e_{l_{k}})$,
$k=1,\ldots ,n$.

\sprth {\bf Proof.}
Let $q\leq j$ and $\phi\in K_{q}$. Assume that
$\phi\in K^{m}\stm K^{m-1}$ for some $m\geq 0$ and let
$\{ K^{s}(\phi )\}_{s=0}^{m}$ be an  analysis of $\phi$.
Let $x_{k}', x_{k}''$ be the initial and final part of $x_{k}$
\wrt $\{ K^{s}(\phi )\}_{s=0}^{m}$.

We shall define $\psi ',\psi ''\in K_{q}$ \st \fa $k$,
$|\phi (x_{k}')|\leq\psi'(m_{j}e_{l_{k}})$ and
$|\phi (x_{k}'')|\leq\psi''(m_{j}e_{l_{k}})$.

\sprth
{\bf Construction of $\psi '$.}

\sprt
For $f\in \cup_{s=0}^{m}K^{s}(\phi )$, we set
\[ D_{f}=\left\{ k: \spp\phi\cap\spp x_{k}'=
 \spp f\cap\spp x_{k}'\right\}.\]
By induction on $s=0,\ldots ,m$, we shall define \fe
$f\in\cup_{s=0}^{m}K^{s}(\phi )$ a function $g_{f}$
with the following properties:

(a) $g_{f}$ is supported on $\{ l_{k}:k\in D_{f}\}$ .

(b) For $k\in D_{f}$, $|f(x_{k}'|\leq m_{j}g_{f}(e_{l_{k}})$

(c) $g_{f}\in K$. Moreover, if $q\leq j$ and $f\in K_{q}$,
then $g_{f}\in K_{q}$.

\sprth
For $s=0$, $f=\pm e^{*}_{m}\in K^{0}(\phi )$,
$D_{f}\neq\emptyset$ only if for some $k$, $l_{k}=m$
and $x_{k}=e_{l_{k}}$. We then set $g_{f}=e^{*}_{l_{k}}$.

Let $s>0$. Suppose that $g_{f}$ have been defined for all
$f\in \cup_{t=0}^{s-1}K_{t}(\phi )$.
Let $f=\frac{1}{m_{q}}(f_{1}+\cdots +f_{d})=
K^{s}(\phi )\stm K^{s-1}(\phi )$, 
where $f_{i}\in K^{s-1}(\phi )$,
$i=1,\ldots ,d$, and  $\{\spp f_{1}<\cdots <\spp f_{d}\}$
is $\calm_{q}$-admissible.
Let $I=\{ i:1\leq i\leq d, D_{f_{i}}\neq\emptyset\}$.
Let $T=D_{f}\stm\cup_{i\in I}D_{f_{i}}$.

Suppose first that $q\leq j$. We set
\[ g_{f}=\frac{1}{m_{q}}\left(\sum_{i\in I}
g_{f_{i}}+\sum_{k\in T}e^{*}_{l_{k}}\right) .\]
Property (a) is obvious. For (b) we have:

\noindent If $k\in D_{f_{i}}$
for some $i\in I$,
\[ |f(x_{k}')|=\frac{1}{m_{q}}|f_{i}(x_{k}')|\leq
\frac{1}{m_{q}}g_{f_{i}}(m_{j}e_{l_{k}})=g_{f}(m_{j}e_{l_{k}}),\]
using the inductive hypothesis. 

\noindent For $k\in T$ we get
\[ |f(x'_{k})|=\frac{1}{m_{q}}\left|\sum f'_{i}(x_{k})\right|
\leq 1\leq
\frac{m_{j}}{m_{q}}=\frac{1}{m_{q}}e_{l_{k}}^{*}(m_{j}e_{l_{k}})
=g_{f}(m_{j}e_{l_{k}}).\]

To show that $g_{f}\in K_{q}$, we need to show that the set
$\{\spp g_{f_{i}} :i\in I\}\cup\{\{l_{k}\}:k\in T\}$ is
$\calm_{q}$-admissible.

Let $G=\{ t_{1}<t_{2}<\cdots <t_{r}\}$ be an ordering of the set
$\{ l_{k}:k\in T\}\cup\{\min\{l_{k}:k\in D_{f_{i}}\},i\in I\}$.
Let $F=\{ s_{1}<s_{2}<\cdots <s_{d}\}\in\calm_{j}$ be
\st $s_{1}\leq\spp f_{1}<s_{2}\leq\spp f_{2}<
\cdots <s_{d}\leq\spp f_{d}$.
By the definition of $x_{k}'$, if $k\in T$ there is
$f_{i}\in\{ 1,\ldots ,d\}\stm I$ \st
$\spp f_{i}\cap\spp x_{k}'\neq\emptyset$,
$\spp f_{i}\cap\spp x_{m}'=\emptyset$ for all $m\neq k$.
This shows that  $r\leq d$ and  $s_{p}\leq t_{p}$
for all $p\leq r$. Then by Lemma 2.1, $G\in\calm_{q}$.

Suppose now that $q>j$. Then we set
$ g_{f}=\frac{1}{m_{j}}\left(\sum_{i\in I}
g_{f_{i}}+\sum_{k\in T}e^{*}_{l_{k}}\right) .$
Since $\{ l_{r},\ldots ,l_{k}\}\in \calm_{j}$, it is 
obvious that $g_{f}\in K$.

Properties (a) and (b) are also easily checked.

The construction of $\psi ''$ is similar.

Finally, we set $\psi =\frac{1}{2}(\psi '+\psi '').$

\sprt
{\bf 2.11 Corollary.} Let $j\in \NN$, $0<\epsilon <
\frac{1}{m_{j}^{3}}$. 
Let $\sum_{k=1}^{n}a_{k}x_{k}$ be an $(\epsilon ,j)$-s.c.c.
Then, for $q<j$, $\phi\in K^{q}$,
$|\phi (\sum a_{k}x_{k})|\leq\frac{4}{m_{q}}$.

\sprt
{\bf Proof.} Combine Propositions 2.5 and 2.10.

\sprt
{\bf 2.12 Definition.} For $j=1,2,\ldots$, $\epsilon >0$,
a finite block \sq $\{ y_{k}\}_{k=1}^{n}$ is said to be an
{\em $(\epsilon ,j)$-rapidly \isq} if the following are satisfied:

\begin{list}{(\alph{lista})} {\usecounter{lista}}
\item There exist $\{ a_{k}\}_{k=1}^{n}$
with $a_{k}\geq 0$, $\sum a_{k}=1$ \st
$\sum_{k=1}^{n}a_{k}y_{k}$ is an $(\epsilon ,j)$-s.c.c.
\item There exist $j_{1},\ldots ,j_{n}$ such that:

 (i) $j+2<2j_{1}<\cdots ,<2j_{n}$, 

(ii) each $y_{k}$ is a normalized
$\left(\frac{1}{m^{4}_{2j_{k}}},2j_{k}\right)$-s.c.c.

(iii) the $\ell^{1}$-norm of $y_{k}$ is dominated by
$\frac{m_{2j_{k+1}}}{m_{2j_{k+1}-1}}$.
\end{list}

The convex combination $y=\sum_{k=1}^{n}a_{k}y_{k}$,
where $\{ a_{k}\}_{k=1}^{n}$ is as above, is said to be an
{\em $(\epsilon ,j)$-rapidly increasing s.c.c.}

\sprt {\bf 2.13 Proposition.}
Let $j\geq 1$. Let $\{ y_{k}\}_{k=1}^{n}$ be an
$(\epsilon ,j)$-rapidly increasing sequence
and $(l_i)_{i=1}^n$ be such that 
$\spp y_1\leq l_1<\spp y_2\leq l_2<\cdots 
\leq l_{n-1}<\spp y_n\leq l_n$ and 
$\{ l_1, \ldots ,l_n\} \in\calm_j$.
Let $j_{k}$ be as in Definition 2.12.
Then, \fe $\phi\in K_{r}$ \txs $\psi\in co(K)$, \st for
$k=1,\ldots ,n$,

$|\phi (y_{k})|\leq 8\psi (e_{l_{k}})$.
Moreover, 

if $r< 2j_{1}$ or $r>2j_{n}$ then
$\psi\in co K_{r}$, 

if $2j_{1}\leq r\leq 2j_{n}$ then
$\psi $ is of the form $\psi=\frac{1}{2}\psi_{1}+
\frac{1}{2}|\phi (y_{k})|e_{l_{k}}$, where $\psi_{1}\in co(K_{r-1})$.

\sprth
{\bf Proof.} The construction is similar to the one
in the proof of Proposition 2.10.

Let $\phi\in K_{r}$. Assume that
$\phi\in K^{m}\stm K^{m-1}$ and let
$\{ K^{s}(\phi )\}_{s=0}^{m}$ be an analysis of $\phi$.
Let $y_{k}'$ and $y_{k}''$ be the initial and final part of $y_{k}$
\wrt $\{ K^{s}(\phi )\}_{s=0}^{m}$.

We shall define $\psi '$ and $\psi ''$ so that 
$|\phi (y_{k}')|\leq 4\psi'(e_{l_{k}})$ and
$|\phi (y_{k}'')|\leq 4\psi''(e_{l_{k}})$.

\sprth
{\bf Construction of $\psi '$}

\sprt
For $f\in \cup_{s=0}^{m}K^{s}(\phi )$, we set
\[ D_{f}=\left\{ k: \spp\phi\cap\spp y_{k}'=
 \spp f\cap\spp y_{k}'\right\}.\]
By induction on $s=0,\ldots ,m$, we shall define \fe
$f\in\cup_{s=0}^{m}K^{s}(\phi )$ a function $g_{f}$
with the following properties:

\begin{list}{\alph{lista})} {\usecounter{lista}}

\item $g_{f}$ is supported on $\{ e_{l_{k}}:k\in D_{f}\}$ 

\item $|f(y_{k}'|\leq 4g_{f}(e_{l_{k}})$ for $k\in D_{f}$.

\item $g_{f}\in K$. Moreover, $g_{f}\in K_{q}$, if $q<2j_{1}$. 

$g_{f}=\frac{1}{2}\psi_{1}+\frac{1}{2}e_{l_{k}}$, with
$\psi_{1}\in K_{s-1}$, $l_{k}\not\in\spp \psi_{2}$.
\end{list}

Let $s>0$. Suppose that $g_{f}$ have been defined for all
$f\in \cup_{t=0}^{s-1}K^{t}(\phi )$.
Let $f=\frac{1}{m_{q}}(f_{1}+\cdots +f_{d})\in
K^{s}(\phi )\stm K^{s-1}(\phi )$, 

{\sc Case 1.} $q<2j_{1}$.

Let $I=\{ i:1\leq i\leq d, D_{f_{i}}\neq\emptyset\}$
and $T=D_{f}\stm\cup_{i\in I}D_{f_{i}}$.
We set
\[ g_{f}=\frac{1}{m_{q}}\left(\sum_{i\in I}
g_{f_{i}}+\sum_{k\in T}e^{*}_{l_{k}}\right) .\]
Properties  (a) and (b) for the case $k\in\cup_{i\in I}D_{f_{i}}$
follow easily from the inductive assumption.
For $k\in T$ we get
\[ |f(y_{k})|=\frac{1}{m_{q}}\left|\sum f_{i}(y_{k})\right|
\leq\frac{4}{m_{q}}\leq 4g_{f}(e_{l_{k}}),\]
by Corollary 2.11, since $q<2j_{k}$ for all $k$.

The proof that $g_{f}\in K_{q}$ is as in the proof
of Proposition 2.10 (Case $q<j$).

{\sc Case 2.} $q\geq 2j_{1}$.

Let $1\leq t\leq n$ be \st $2j_{t}\leq q<2j_{t+1}$.
We shall define $g_{f}',g_{f}'', g_{f}'''$ supported on
$\{ l_{k}:k\leq t-1\}$, $\{ l_{t}$\} and
$\{ l_{k}:k\geq t+1\}$ respectively.
Let \[ f'=f|\bigcup_{k=1}^{t-1}\spp y_{k}',\;\;
f''=f|\spp y_{t}',\;\;
f'''=f|\bigcup_{k=t+1}^{n}\spp y_{k}'.\]
Let \[ I'=\left\{\rule{0mm}{4mm} i:1\leq i\leq d, \{ k\leq t-1\}
\cap D_{f_{i}}\neq\emptyset\}\right\} \;\;{\rm and}\;\;
T'=D_{f'}\stm\bigcup_{i\in I'}D_{f_{i}}\]
and define similarly $I'''$ and $T'''$. We set:

(i)\hfill${\displaystyle g_{f}'=\frac{1}{m_{q-1}}
\left(\sum_{i\in I'}g_{f_{i}}+\sum_{k\in K'}
e^{*}_{l_{k}}\right).}$\hfill\mbox{ }\\
It is obvious that $g'_{f}\in K_{q-1}$. Moreover,
for $k\leq t-1$ we get:\\
If $k\in\cup_{i\in I'}D_{f_{i}}$, then 
$|f'(y'_{k})|\leq 4g_{f}(e_{l_{k}})$ by the inductive assumption.\\
If $k\in T'$, then
\beq |f'(y_{k}')|&=&\frac{1}{m_{q}}\left|\left(
\sum f_{i}\right) (y'_{k})\right|\leq\frac{1}{m_{q}}\| y_{k}\|_{1}\\
&\leq&\frac{m_{2j_{t}}}{m_{2j_{t}-1}}\cdot\frac{1}{m_{q}}
\leq\frac{m_{q}}{m_{q-1}}\cdot\frac{1}{m_{q}}=
\frac{1}{m_{q-1}}=g'_{t}(e_{l_{k}}).\eeq

(ii)\hfill${\displaystyle g_{f}''=|f(y_{t})|\cdot e^{*}_{l_{t}}
.}$\hfill\mbox{}\\

(iii)\hfill${\displaystyle g_{f}'''=\frac{1}{m_{q-1}}
\left(\sum_{i\in I'''}g_{f_{i}}+\sum_{k\in K'''}
e^{*}_{l_{k}}\right).}$\hfill\mbox{}\\
Then $g_{f}'''\in K_{q-1}$.
For $k\in\cup_{i\in I'''}D_{f_{i}}$,  
$|f'''(y'_{k})|\leq 4g_{f}(e_{l_{k}})$ by the inductive assumption.
For $k\in T'''$ we get
\[ |f'''(y_{k}')|=\frac{1}{m_{q}}\left|
\sum (f_{i}) (y'_{k})\right|\leq\frac{4}{m_{q}}\leq\frac
{1}{m_{q-1}}=g'''_{f}(e_{l_{k}}),\]
using Corollary 2.11. Finally, we set
\[ g_{f}=\frac{1}{2}\left( g_{f}'+g_{f}''+g_{f}'''\right) .\]
This completes the proof for $\psi '$.

The construction of $\psi ''$ is similar.

Finally, $\psi =\frac{1}{2}(\psi '+\psi '')$.

\sprt
{\bf 2.14 Proposition.} Let $\sum_{k=1}^{n}a_{k}y_{k}$ be a
$\left(\frac{1}{m^{4}_{j}},j\right)$-rapidly increasing s.c.c.
Then for $i=0,1,2,\ldots$, $\phi$ in $K_{i}$, we
have the following estimates:

\begin{list}{(\alph{lista})} {\usecounter{lista}}
\item $|\phi (\sum_{k=1}^{n}a_{k}y_{k})|\leq
\frac{16}{m_{i}m_{j}}$, if $i<j$,
\item $|\phi (\sum_{k=1}^{n}a_{k}y_{k})|\leq
\frac{8}{m_{i}}$, if $j\leq i<2j_{1}$ or $2j_{n}<i$,
\item $|\phi (\sum_{k=1}^{n}a_{k}y_{k})|\leq
\frac{4}{m_{i-1}}+4|a_{k_{0}}||\phi(y_{k_{0}})|$, if 
$2j_{k_{0}}\leq i<2j_{k_{0}+1}$.
\end{list}

\sprt{\bf Proof.} It follows easily from Proposition 2.5 and
Lemma 2.13.

\sprt {\bf 2.15 Corollary.} If $\sum_{k=1}^{n}a_{k}y_{k}$ is a 
$\left(\frac{1}{m^{4}{j}},j\right)$-rapidly increasing s.c.c.
then

\mbox{ }\hfill $\displaystyle
\frac{1}{4m_{j}}\leq\left\|\sum_{k=1}^{n} a_{k}y_{k}\right\|\leq
\frac{8}{m_{j}}$\hfill\mbox{ }.

\sprt{\bf 2.16 Corollary.} $X_{u}$ is arbitrarily distortable.

\newcommand{\nrm}{|\! |\! |}
\sprth{\bf Proof.} Choose $i_{0}$ arbitrarily large.
Let \[ \nrm x\nrm =\frac{1}{m_{i_{0}}}\| x\| +
\sup\left\{ \phi (x):\phi\in K_{i_{0}}\right\} .\]
Let $Y$ be a block subspace of $X_{u}$. Let $j>i_{0}$.
Using Lemma 2.9, we can choose the following vectors in $Y$,
\beq y&=&\sum_{k=1}^{n}a_{k}y_{k},\;
{\rm  a}\; 
\left(\frac{1}{m^{4}_{j}},
j\right)\mda{\rm rapidly\; increasing\; s.c.c.}\\
z&=&\sum_{l=1}^{m}b_{l}z_{l},\;
{\rm  a}\; 
\left(\frac{1}{m^{4}_{i_{0}}},i_{0}\right)\mda
{\rm rapidly\; increasing\; s.c.c.}
\eeq
Then, by Corollary 2.14,
\[ \nrm m_{j}y\nrm\leq\frac{8}{m_{i_{0}}}+
\frac{16}{m_{i_{0}}}=\frac{24}{m_{i_{0}}}
\;{\rm while}\; \|m_{j}y\|\geq\frac{1}{4},\]
\[ \nrm m_{i_{0}}z\nrm\geq\frac{1}{4}\;{\rm while}\;
\| m_{i_{0}}z\|\leq 8.\]
This completes the proof.

\section{The space $X$}

We turn now to defining the Banach space $X$ not containing 
any unconditional basic sequence.
The norm of the space is related to that of $X_{u}$ 
introduced in the previous section.
Specifically, the norm will be defined by a family
$\{ L_{j}\}_{j=0}^{\infty}$ of subsets of $c_{00}$ \st
each $L_{j}$ is contained in the set $K_{j}$ used in
the definition of $X_{u}$. 
We consider the countable set
\[ G=\left\{ \left( x^{*}_{1},x^{*}_{2},\ldots ,x^{*}_{k}
\right) :x^{*}_{i}\in\bigcup_{j=0}^{\infty}K_{j},
x^{*}_{1}<x^{*}_{2}<\cdots <x^{*}_{k}\right\} .\]
There exists a function $\Phi :G \lra\{ 2j\}_{j=0}^{\infty}$
one to one \st if $\left( x^{*}_{1},\ldots ,x^{*}_{k}
\right)\!\in G$, $x^{*}_{1}\in K_{j_{1}},\ldots ,
x^{*}_{k}\in K_{j_{k}}$, then
\[ \Phi\left( (x^{*}_{1},\ldots ,x^{*}_{k})
\right) >\max\{ j_{1},\ldots ,j_{k}\}.\]
For $n=0,1,2,\ldots$ we define by induction sets
 $\{ L^{n}_{j}\}_{j=0}^{\infty}$ \st $L^{n}_{j}$
is a subset of $K^{n}_{j}$ and
$\{ L^{n}_{j}\}_{n=0}^{\infty}$ is an increasing family.
We set $L^{0}_{j}=\{\pm e_{n}:n=1,2,\ldots \}$.
Suppose that $\{ L^{n}_{j}\}_{j=0}^{\infty}$ have been defined
and set
\beq 
L^{n+1}_{2j}&=&\pm L^{n}_{2j}\cup\left\{\frac{1}{m_{2j}}
\left( x^{*}_{1}+\cdots +x^{*}_{d}\right) :
x_{i}\in\bigcup_{j=0}^{\infty}L^{n}_{j},
\left(\spp x^{*}_{1},\ldots ,\spp x^{*}_{d}\right)
\right.\\
& &\;\;\;\;\;\;\;\;\;\;{\rm is}
\left. \calm_{2j}\; {\rm \mda admissible}\rule{0mm}{7mm}\right\} ,
\\
L'^{n+1}_{2j+1}&=&\pm L^{n}_{2j+1}\cup\left\{\frac{1}{m_{2j+1}}
\left( x^{*}_{1}+\cdots +x^{*}_{d}\right) :
x_{1}^{*}\in L^{n}_{2k}, k>2j+1,\right.\\ 
& & \;\;\;\;\;\;\;\;x^{*}_{i}\in 
L^{n}_{\Phi (x^{*}_{1},\ldots ,x^{*}_{i-1})}
\; {\rm for}\; 1<i\leq d\;{\rm and}\; 
\left(\spp x^{*}_{1},\ldots ,\spp x^{*}_{d}\right)\\
& &\left.\;\;\;\;\;\;\;\;{\rm is\; 
\calm_{2j+1}\mda admissible}\rule{0mm}{6mm}\right\}\;{\rm and}\\
L^{n+1}_{2j+1}&=&\left\{ E_{s}x^{*} :x^{*}\in L'^{n+1}_{2j+1},
s \in \NN , E_{s} =\{ s, s+1, \ldots\}\right\} 
\eeq
This completes the definition of $L^{n}_{j}$,
$n=0,1,2,\ldots$, $j=0,1,2,\ldots$.
It is obvious that each $L^{n}_{j}$ is a subset of 
the corresponding set $K^{n}_{j}$.

We set $L_{j}=\cup_{n=0}^{\infty}L^{n}_{j}$ and we
consider the norm on $c_{00}$ defined by the family
$L=\cup_{j=0}^{\infty}L_{j}$.
The space $X$ is the completion of $c_{00}$ under this norm.
It is easy to see that $\{ e_{n}\}_{n=1}^{\infty}$
is a bimonotone basis of $X$.

\sprt {\bf 3.1 Remark.}
An alternative implicit definition of the norm of the space
$X$ is the following. For $x\in c_{00}$,
\beq \| x\|\! =\! \max\left\{\| x\|_{0},
\sup\left\{\frac{1}{m_{2j}}\sum_{k=1}^{n}\| E_{k}x\| ,\rule{0mm}{7mm}
j\in\NN, n\in\NN ,\{ E_{1}<\cdots <E_{n}\}
\;{\rm is}\right.\right.\\
\left.\left.\;\; \calm_{2j}{\rm \mda admissible}\right\} ,
\sup\left\{ |\phi (x)|:\phi\in\bigcup_{j=1}^{\infty}L^{2j+1}\right\}
\right\} .\eeq
Hence, for $j=1,2,\ldots$ 
and for
$x_{1}<x_{2}<\cdots ,<x_{n}$ in $c_{00}$ \st
$\{ \spp x_{1},\spp x_{2},\ldots ,\spp x_{n}\}$ is
$\calm_{2j}$-admissible, we have that
$\|\sum_{k=1}^{n}x_{k}\|\geq\frac{1}{m_{2j}}
\sum_{k=1}^{n}\| x_{k}\|$.
This allows us to have the following result in the same manner
as Lemma 2.9.

\sprt
{\bf 3.2 Lemma.}
For $j=1,2,\ldots$ and every normalized block \sq
$\{ x_{k}\}_{k=1}^{\infty}$ in $X$ \txs $\{ y_{s}\}_{s=1}^{n}$
a finite block \sq of $\{ x_{k}\}_{k=1}^{\infty}$ \st
$\sum_{s=1}^{n}a_{s}y_{s}$ is a normalized
$\left(\frac{1}{m^{4}_{j}},2j\right)$-s.c.c.

\sprt {\bf 3.3 Proposition.}
Let $\sum a_{k}x_{k}$ be an 
$\left(\frac{1}{m^{3}_{j}},j\right)$-s.c.c.
defined by an $\left(\frac{1}{m^{3}_{j}},j\right)$-basic s.c.c.
$\sum_{k=1}^{n}a_{k}e_{l_{k}}$. Then \fe $s\leq j$ and
$\phi$ in $L_{s}$ \txs $\psi$ in $K_{s}$ \st
\[ \phi\left(\sum_{k=1}^{n}a_{k}x_{k}\right)\leq
2\psi\left(\sum_{k=1}^{n}a_{k}m_{j}e_{l_{k}}\right).\]
The proof of this is similar to the proof
of Proposition 2.10. 

\sprt
{\bf 3.4 Proposition.} Let $\sum_{k=1}^{n}a_{k}y_{k}$ be a
$\left(\frac{1}{m^{4}{j}},j\right)$-rapidly increasing s.c.c.
in $X$. Then for $i=I$, $\phi\in L_{i}$, we
have the following estimates:

\begin{list}{(\alph{lista})} {\usecounter{lista}}
\item $|\phi (\sum_{n=1}^{k}a_{k}y_{k})|\leq
\frac{16}{m_{i}m_{j}}$ if $i<j$,
\item $|\phi (\sum_{n=1}^{k}a_{k}y_{k})|\leq
\frac{8}{m_{i}}$ if $j\leq i<2j_{1}$ or if $i>2j_{n}$,
\item $|\phi (\sum_{n=1}^{k}a_{k}y_{k})|\leq
\frac{4}{m_{i-1}}+4|a_{k_{0}}||\phi(y_{k_{0}})|$ if 
$2j_{k_{0}}\leq i<2j_{k_{0}+1}$.
\end{list}

In particular, $\|\sum_{k=1}^{n}a_{k}y_{k}\|\leq\frac{8}{m_{j}}$.

\sprth
This is proved similarly to Proposition 2.14.

\sprt{\bf 3.5 Lemma.}
 Let $\sum_{k=1}^{n}a_{k}x_{k}$ be an
$(\epsilon ,j)$-s.c.c. and $2i<j$. Then \fe $\psi$ in
$K_{i}$, $\psi =\frac{1}{m_{2i}}\left(x^{*}_{1}+
\cdots +x^{*}_{d}\right)$ and every choice
$\{ a_{k_{t}}\}_{t=1}^{s}$ \st \fe
$t=1,\ldots ,s$ \tx $\tau_{1}<\tau_{2}\leq d$ so that
$x^{*}_{\tau_{1}}(x_{k_{t}})\neq 0$ and
$x^{*}_{\tau_{2}}(x_{k_{t+1}})\neq 0$, we have
$\sum_{t=1}^{s}a_{k_{t}}<m_{2i}\cdot\epsilon$.

\sprth{\bf Proof.} Let $\sum_{k=1}^{n}a_{k}e_{l_{k}}$
be the $(\epsilon ,j)$-basic s.c.c. that defines the vector
$\sum_{k=1}^{n}a_{k}x_{k}$.
It is easy to check that the set
$(e_{l_{k_{t}}})_{t=1}^{s}$ is $\calm_{2i}$-admissible, hence
\[ \frac{1}{m_{2i}}\sum_{t=1}^{s}a_{k_{t}}\leq
\left\|\sum_{t=1}^{s}a_{k_{t}}e_{l_{k_{t}}}\right\|_{2i}
<\epsilon\] 
and the proof is complete.

\sprt{\bf 3.6 Proposition.}
Let $j>100$ and $\{ j_{1},\ldots ,j_{n}\}$ be \st
$2j+1<j_{1}<j_{2}<\cdots <j_{n}$ in {\bf N}.
Suppose that $\{ y_{k}\}_{k=1}^{n}$,
 $\{ \theta_{k}\}_{k=1}^{n}$ are such that:

\begin{list}{(\roman{lista})} {\usecounter{lista}}
\item Each $y_{k}$ is a
$\left(\frac{1}{m^{4}_{2j_{k}}},2j_{k}\right)$-rapidly increasing s.c.c.
\item $y^{*}_{k}\in L_{2j_{k}}$,
$y^{*}_{k}(y_{k})\geq\frac{1}{2m_{2j_{k}}},$
$y^{*}_{k_{1}}(y_{k_{2}})=0$ for $k_{1}\neq k_{2}$,
$y^{*}_{1}<y^{*}_{2}<\cdots <y^{*}_{n}$.
\item $\frac{1}{8}\leq\theta_{k}\leq 2$ and
$y^{*}_{k}(m_{2j_{k}}\theta_{k}y_{k})=1$.
\item $2j_{k}=\Phi \left( y^{*}_{1},\ldots ,y^{*}_{k-1}\right)$,
$k=2,\ldots ,n$.
\item There exists decreasing \sq $\{ a_{k}\}_{k=1}^{n}$ \st
$\sum_{k=1}^{n}a_{k}y_{k}$ is\linebreak
$\left(\frac{1}{m^{4}_{2j+2}},2j+1\right)$-s.c.c.
\end{list}

Let $(\epsilon_{k})_{k=1}^{n}$ be \st $\epsilon_{k}=1$
if $k$ is even and $\epsilon_{k}=-1$ if $k$ is odd.
It is clear that
\[ \left|\sum_{k=k_{1}}^{k_{2}}\epsilon_{k}a_{k}m_{2j_{k}}\theta_{k}
y^{*}_{k}(y_{k})\right|\leq a_{k_{1}}\leq\frac{1}{m^{4}_{2j+2}}.\]

Then $\|\sum\epsilon_{k}a_{k}m_{2j_{k}}y_{k}\|
\leq \frac{18}{m_{2j+2}}$.

\sprth
The proof is given in several steps. Our aim is to 
show that \fe 

$\phi\in\cup_{j=0}^{\infty}L_{j}$,
\[\phi \left(\sum\epsilon_{k}a_{k}m_{2j_{k}}\theta_{k}y_{k}\right)
\leq \frac{18}{m_{2j+2}}.\]

{\bf 3.7 Lemma.} Let $y$ be a
$\left(\frac{1}{m^{4}_{2j}},2j\right)$-rapidly increasing s.c.c.
and $z_{1}^{*},\ldots ,z_{d}^{*}$ be in 
$L_{2j_{1}},\ldots ,L_{2j_{2}}$, respectively, \st
$2j_{k}\neq 2j$ for all $k=1,\ldots ,d$ and
$(\spp z_{1},\ldots ,\spp z_{d})$ is $\calm_{i}$-admissible
for some $i<\min\{ 2j,2j_{i},\ldots ,2j_{d}\}$.
Then \[ \left| z_{1}^{*}+\cdots +z_{d}^{*}(y)\right|
<\sum_{k=1}^{d_{1}}\frac{16}{m_{2j_{k}}m_{2j}}+
\sum_{k=d_{1}+1}^{d}\frac{8}{m_{2j_{k}-1}}+
\frac{1}{m^{2}_{2j}},\]
where $j_{1}<\cdots <j_{d_{1}}<j<j_{d_{1}+1}<\cdots <j_{d}$.

\sprth {\bf Proof.}
Let $y=\sum_{n=1}^{l}a_{n}x_{n}$ be the expression of $y$
as a rapidly increasing 
$\left(\frac{1}{m^{4}_{2j}},2j\right)$-s.c.c.
First we notice that for $1\leq k<d_{1}$,
$|z_{k}^{*}(y)|\leq\frac{16}{m_{2j_{k}}m_{2j}}$
and \fe $d_{1}+1\leq k\leq d$,
$|z_{k}^{*}(y)|\leq\frac{4}{m_{2j_{k}-1}}+
4a_{n_{k}}|z_{k}^{*}(x_{n_{k}})|$.

For $n=1,2,\ldots ,l$ we set
$I_{n}=\{ k:n_{k}=n\}$.
Then $|\sum_{k\in I_{n_{0}}}z_{k}^{*}(y)|\leq 
\sum_{k\in I_{n_{0}}}\frac{4}{m_{2j_{k}-1}}+
64a_{n_{0}}$ since
 $|\sum_{n\in I_{n_{0}}}z_{k}^{*}(y)|\leq 16$.
Therefore \beq
\left|\sum_{k=d_{1}+1}^{d}z_{k}^{*}(y)\right|&\leq&
\sum_{n=1}^{l}\left|\sum_{k\in I_{n}}z_{k}^{*}(y)
\right| +
\left|\sum_{k\not\in\cup_{n=1}^{l}I_{n}}z_{k}^{*}(y)
\right| \\
&\leq&\sum_{k=d_{1}+1}^{d}\frac{4}{m_{2j-1}}+
64\sum_{k\in S}a_{n_{k}},
\eeq
where $S\sbs\{ d_{1}+1,\ldots ,d\}$ is \st for
$k_{1}\neq k_{2}$ in $S$, $n_{k_{1}}\neq
n_{k_{2}}$ and $S$ is maximal with this property.
Observe that $\{\spp x_{n_{k}}\}_{k\in S}$ is
$\calm_{i}$-admissible. Hence, by Lemma 3.5,
we get that $16\sum_{k\in S}a_{n_{k}}<
16m_{i}\cdot\frac{1}{m^{4}_{2j}}\leq
\frac{1}{m^{2}_{2j}}$. 
This completes the proof of the lemma.

\sprt{\bf 3.8 Proposition.}
Let $y_{k}$ be a
$\left(\frac{1}{m^{4}_{2j_{k}}},2j_{k}\right)$-rapidly increasing s.c.c.
for $k=1,\ldots ,n$ and$z_{s}^{*}$ be in 
$L_{2t_{s}}$ for  $s=1,\ldots ,d$.
Suppose that $2j+1<j_{1}<\cdots<j_{n}$,
$2j+1<t_{1}<\cdots <t_{d}$,
$\{ j_{1},\ldots ,j_{n}\}\cap \{ t_{1},\ldots ,t_{d}\}
=\emptyset$ and
$(\spp z_{1}^{*},\ldots ,\spp z_{d}^{*})$ is $\calm_{2j+1}$-admissible.\\
Then ${\displaystyle \left(\sum_{s=1}^{d}z_{s}^{*}\right)
\left(\sum_{k=1}^{n}a_{k}m_{2j_{k}}y_{k}\right)
<\sum_{k=1}^{n}|a_{k}|\frac{1}{m^{2}_{2j+2}} }$.

\sprt{\bf Proof.} It follows easily from Lemma 3.7.

\sprt
{\bf 3.9 Proposition.} For every $\phi$ in
$L_{2j+1}$ we have
\[\phi \left(\sum_{k=1}^{n}
\epsilon_{k}a_{k}m_{2j_{k}}\theta_{k}y_{k}\right)
\leq \frac{4}{m^{2}_{2j+2}}.\]

\sprth {\bf Proof.} Let $\phi =\frac{1}{m_{2j+1}}
(x_{k_{1}-1}^{*}+y^{*}_{k_{1}}\cdots +y_{k_{2}}^{*}
+z_{k_{2}+1}^{*}+\cdots
+z_{d}^{*})$ be in $L_{2j+1}$. 
Set $w=\sum_{k=1}^{n}\epsilon_{k}a_{k}m_{2j_{k}}\theta_{k}y_{k}$.
Then \beq |\phi (w)|\!\! &\leq&\! \frac{1}{m_{2j+1}}
|(y_{k_{1}}^{*}+\cdots +y_{k_{2}}^{*})(w)|+|x_{k_{1}-1}^{*}(w)|+
|(z_{k_{2}+2}^{*}+\cdots +z_{d}^{*})(w)|\\
\!\!&<&\! a_{k_{1}}+a_{k_{1}-1}+\frac{1}{m^{2}_{2j+2}}
+\frac{2}{m^{2}_{2j+2}}<\frac{4}{m^{2}_{2j+2}}.
\eeq

\sprt{\bf 3.10 Proposition.}
If $i>2j+1$ and $\phi$ is in $L_{i}$,
we have the following
\[ \phi \left(\sum_{k=1}^{n}a_{k}m_{2j_{k}}y_{k}\right)
\leq\left\{\begin{array}{ll}
\sum_{k=1}^{n}|a_{k}|\frac{16}{m_{i}}\;\;
{\rm if}\; 2j+1<i<2j_{1}\;{\rm or}\; i>2j_{n}\arrsp
4\sum_{k=1}^{n}\frac{|a_{k}|}{m_{i-2}}+
8|a_{k_{0}}| \;\;{\rm if}\; 2j_{k_{0}}\leq
i< 2j_{k_{0}+1}\end{array}\right. .\]

\sprth
{\bf Proof} The case $2j+1<i<2j_{1}$ follows from 
Proposition 3.4 (a),
the case 
$2j_{k_{0}}\leq i\leq 2j_{k_{0}+1}$ 
and $i>2j_{n}$ from Proposition 3.4 (b)
and (c) and the fact that for $i_{1}<i_{2}$,
$\frac{m_{i_{1}}}{m_{i_{2}}}<\frac{1}{m_{i_{2}-1}}$.

\sprt{\bf 3.11 Lemma.}
For every $i<2j+1$ and $\phi $ in $L_{i}$ 
we have 
\[ \phi \left(\sum \epsilon _{k}a_{k}m_{2j_{k}}\theta _{k}y_{k}
\right)\leq \frac{18}{m_{2j+2}}.\]

{\bf Proof.}
Let $\{ K^{s}(\phi )\}_{s=0}^{m}$ be an analysis of the functional
$\phi $. We set

$ W=\left\{ \rule{0mm}{4mm}k:\exists s>s_{k}{\rm\; and}\;
f^{k}\in K^{s}(\phi )\; {\rm such\; that}\; f=\right.$

$\left. 
\;\;\;\;\frac{1}{2j+1}\left(x^{*}_{t_{1}-1}+y^{*}_{t_{1}}+\cdots
+y^{*}_{t_{2}}+z ^{*}_{t_{2}+1}+\cdots +z ^{*}_{d}\right)
{\rm \; and}\; t_{1}\leq k\leq t_{2}\right\}. $

\noindent For every $k$ in $W$ we denote by 
$f^{k}$ the function $f$ which witnesses the belongingness
of $k$ in $W$ and is \st if $f^{k}$ is in $K^{s}(\phi )$,
$s$ is the maximum level where such an $f$ occurs.
We also  denote by $t_{1}^{k}, t^{k}_{2}$ the corresponding
$t_{1}, t_{2}$ of the function $f^{k}$.
Observe that the family 
$\{ [ t^{k}_{1},t^{k}_{2}]_{k\in W}\}$ defines a partition
of the set $W$.
Therefore we write the above segments as
$\{ T_{\sigma}\}_{\sigma =1}^{w}$.

We set
\beq S_{1}&=&\left\{\rule{0mm}{6mm} k:k\in W^{c}\;{\rm and\; for\; every\;}
f\in \bigcup_{s=0}^{m}K^{s}(\phi )\;{\rm with}\right.\\
& &
\left.\spp \phi \cap \spp y_{k}=
\spp f\cap\spp y_{k}\neq\emptyset\;
{\rm we \; have}\; f\in\bigcup_{i'\leq 2j}L_{i'}\right\}\\
S_{2}&=&\left\{ k: k\in W \;{\rm for\; every}\; f\;{\rm strictly\;
extending}\; f^{k}, f \;{\rm is\; in}\; 
\cup_{i'\leq 2j}L_{i'}\right\}\eeq

We set $I_{1}=S_{1}\cup S_{2}$.

\sprtq {\sc Claim 1:} There exists $\psi$ with 
$\|\psi\|^{*}_{2j}\leq 1$ \st

\beq \left|\phi \left(\sum _{k\in I_{1}}
\epsilon _{k}a_{k}m_{2j_{k}}\theta_{k}y_{k}\right)\right|
&\leq& \psi \left(\sum_{k\in S_{1}}a_{k}e_{l_{k}}+
2\sum_{\sigma =1}^{w}a_{k_{\sigma}}
a_{k_{\sigma}}e_{l_{k_{\sigma}}}\right)\\
&\leq& 64\| a_{k}e_{l_{k}}\|_{2j} <\frac{1}{m^{2}_{2j+2}},\eeq 
where 
$a_{k_{\sigma}}=\max\{ a_{k}:k\in T_{\sigma}\}$.

\sprth
{\bf Proof of the claim.} We follow the same
 procedure as in the proof of Proposition 2.10.
For $f$ in $K^{s}(\phi )$ we construct by induction a functional
$g_{f}$ with $\| g_{f}\| ^{*}_{2j}\leq 1$ and \st

If $D_{f}=\{ k:k\in S_{1}, s>s_{k}\}\cup
\{ k_{\sigma}: f\;{\rm extends}\; f^{k_{\sigma}}\}$.

Then 
\begin{list}{\arabic{lista})} {\usecounter{lista}}
\item $g_{f}=0$ if $D_{f}=\emptyset$
\item $\spp g_{f}\sbs \{ l_{k}:k\in D_{f}\}$
\item $|f(a_{k}m_{2j_{k}}\theta_{k}y'_{k})|\leq
16 g_{f}(a_{k}e_{l_{k}})$ for $k\in S_{1}\cap
D_{f}$ ($y_{k}'$ is the innitial part of $y_{k}$), 

$|f(\sum_{k\in T_{\sigma}}\epsilon _{k}a_{k}
m_{2j_{k}}\theta_{k}y_{k})|\leq
2 g_{f}(a_{k_{\sigma}}e_{l_{k_{\sigma}}})$ if
$k_{\sigma}\in D_{f}$, $k=k_{\sigma}$. 
\end{list}

Suppose that $g_{f}$ has been constructed for $f$
in $K^{s'}(\phi )$ for all $s'<s$. Choose $f$ in $K^{s}(\phi )
\stm K^{s'}(\phi )$ with $D_{f}\neq\emptyset$.
Then either $i'\leq 2j$ or $i'=2j+1$.

If $i'\leq 2j$, we observe that 
$D_{f}=\cup_{t=1}^{d}D_{f_{i}}\cup (D_{f}\cap S_{1})$,
where \linebreak $f=\frac{1}{m_{i'}}(f_{1}+\cdots +f_{d})$ and
$g_{f}=\frac{1}{m_{i'}} (\sum_{t\in I} g_{f_{t}}+
\sum_{k\in T}e^{*}_{l_{k}})$
with $I$ and $T$ are as in Proposition 2.10. 
Therefore we conclude  (as in Proposition 2.10) that
$\| g_{f}\| ^{*}_{2j}\leq 1$. Now using the fact that
$i'\leq 2j<2j_{k}$ we conclude that  for $k\in T$
we get the inductive condition (3) for the function $g_{f}$.

If $i=2j+1$ then $f=f^{k_{\sigma}}$, hence $D_{f}=
\{ k_{\sigma}\}$. Set $g_{f}=e ^{*}_{l_{k_{\sigma}}}$ and
$|f(\sum_{k\in T_{\sigma}}\epsilon _{k}a_{k}
m_{2j_{k}}\theta_{k}y_{k})| =
|f\sum_{k\in T_{\sigma}}\epsilon _{k}a_{k}
m_{2j_{k}}y ^{*}_{k}\theta_{k}(y_{k})|
\leq a_{k_{\sigma}}$.

The same proof works for the final parts  
$\{ y_{k}''\}_{k\in S_{1}}$ and the proof of Claim 1 is complete.

\sprtq
We set 
\beq
I_{2}&=&\left\{\rule{0mm}{6mm} k\not
\in W\cup I_{1}\;{\rm and\; there\; exists}\;
s>s_{k}, f\; {\rm in}\; K^{s}(\phi ),\right.\\
& &\left.\!\!\!\!\!\!\!\!\!\spp f\cap \spp y_{k}=\spp \phi \cap \spp y_{k}\;
{\rm such\; that}\; 
| f(m_{2j_{k}}\theta_{k}y_{k})|\leq \frac{16}{m_{2j+2}}\right\} .\\
I_{3}& =&\left\{ k\in W:\;{\rm there\; exists}\;
f\in \cup_{s=1}^{m}K^{s}(\phi ), f\; {\rm strictly}\;
{\rm extending}\; f^{k}\; {\rm and}\right.\\
& &\left. f\in L_{i} \;{\rm for \; some}\; i\geq 2j+1
\rule{0mm}{4mm}\right\}
\eeq
It is  obvious that 
$|f(\sum_{k\in I_{2}}\epsilon _{k}a_{k}
m_{2j_{k}}\theta_{k}y_{k})| <
\sum_{k\in I_2}\frac{2|a_{k}|}{m_{2j+2}}.$

Suppose now that $k\in I_{3}$ and $i=2j+1$. 
Then $f=\frac{1}{m_{2j+1}}(f_{1}+\cdots +f_{d})$ and \txs
$f_{t}$ with $f_{t}$ extending $f^{k_{\sigma }}$ and $f_{t}$
in $L_{2s_{t}}$ for some $s_{t}>2j+1$.
Hence for $k\in T_{\sigma }$,
\[ \left| f(a_{k}m_{2j_{k}}\theta_{k}y_{k})\right|
<\frac{|f^{k_{\sigma }}(a_{k}m_{2j}\theta_{k}y_{k})|}
{m_{2s_{t}}}<\frac{2|a_{k}|}{m_{2j+2}}.\]

Similarly, if $f$ extends $f^{k_{\sigma }}$ and $f$ is in 
$L_{i}$ for $ i\geq 2j+2$,  we get that
for $k\in T_{\sigma }$,
\[ \left| f(a_{k}m_{2j_{k}}\theta_{k}y_{k})\right|
\leq\frac{|f^{k_{\sigma }}(a_{k}m_{2j_{k}}\theta_{k}y_{k})|}
{m_{i}}\leq\frac{|a_{k}|}{m_{2j+2}}.\] 
Hence 
\[  \phi \left( \sum_{k\in I_{3}}\epsilon _{k}
a_{k}m_{2j_{k}}y_{k}\right)
\leq\frac{\sum_{k\in I_{3}}|a_{k}|} {m_{2j+2}}.\]

Notice that $W\sbs I_{1}\cup I_{3}$. We set 
$I_{4}=\{ 1,2,\ldots ,n\}\stm I_{1}\cup I_{2}\cup I_{3}$.

It remains to estimate the quantity
$\phi\left(\sum_{k\in I_{4}}\epsilon_{k}a_{k}m_{2j_{k}} 
\theta_{k}y_{k}\right) .$
To this end, we first make the following observations.
\begin{list}{(\roman{lista})} {\usecounter{lista}}
\item If $k_{0}\in I_{4}$, \txs $f\in\cup_{s>s_{k}}
K_{s}(\phi )$ \st
 
$\spp f\cap\spp y_{k_{0}}=\spp \phi\cap\spp y_{k_{0}}$
$f\in L_{i'}$ with $2j_{k_{0}}\leq i'<2j_{k_{0}+1}$
and 

$f_{i'}(\sum \epsilon _{k}a_{k}m_{2j_{k}}\theta_{k}
y_{k})|\leq
4\sum_{k=1}^{n}\frac{|a_{k}|}{m_{i'-2}}+4|a_{k_{0}}|$
(Proposition 3.10).\\
Otherwise $k$ is in $I_{1}\cup I_{2}\cup I_{3}$. 
\item If $k\in I_{4}$ and $f$ is as above, then \fe
$s'>s$, $f'\in F_{s'}(\phi )$ \st $\spp f\sbs\spp f'$
we have $f'\in\cup_{i'\leq 2j}L_{i'}$. 
Otherwise $k$ is in $I_{2}$.
\item If $k_{1},k_{2}$ are in $I_{4}$, $k_{1}\neq k_{2}$ and
$f\in K^{s}(\phi )$ is such that $s>\max\{ s_{k_{1}},s_{k_{2}}\}$,
 $f\in L_{i'}$ for $i'\geq 2j+1$ 
and 
$f(y_{k_{1}})\neq 0$, then
$f(y_{k_{2}})=0$.
This happens since any such $f$ is pathological 
on only one $y_{k}$, therefore, if $f$ is different than
zero on both, then at least one does not belong to $I_{4}$.
\end{list}

After this we prove the following

\sprth{\sc Claim 2:}
${\displaystyle \left|\phi\left(\sum_{k\in I_{4}}
\epsilon_{k}a_{k}m_{2j_{k}}
\theta_{k}y_{k}\right)\right|\leq 16\left\|\sum_{k\in I_{3}}
a_{k}e_{l_{k}}\right\|_{2j} }.$

\sprth
{\bf Proof:} For $f\in K_{s}(\phi )$ we define, as in Claim 1,
 the set $D_{f}$ and inductively a functional 
$g_{f}$ in $\cup_{i'\leq 2j}K_{i'}$.
If $f$ is in $K_{s}(\phi )$ and $f\in L_{i'}$ for some
$i'\leq 2j$ then we define $g_{f}$ as in Claim 1.
If $i'\geq 2j+1$ then either $D_{f}=\emptyset$ or
$D_{f}=\{ k\}$. This follows immediately from property
(iii) of $I_{4}$ mentioned above.
Therefore, setting $g_{f}=\frac{1}{m_{i'}}e_{k}^{*}$,
we easily check that the inductive assumptions are 
satisfied. The proof of the claim is complete.

\sprth
Hence\hfill$\displaystyle\left|\phi\left(\sum_{k\in I_{4}}
\epsilon_{k}a_{k}m_{2j_{k}}
\theta_{k}y_{k}\right)\right| <\frac{1}{m^{2}_{2j+2}}$\hfill\mbox{ }

The proof is complete.

\sprt {\bf 3.12 Proposition.} Let $(x_{k})_{k\in\NN}$,
$(w_{k})_{k\in{\bf N}}$ be two normalized block
sequences in the space $X$. Then \tx 
$\{ y_{n}\}_{n=1}^{d}$,
$\{ y_{n}^{*}\}_{n=1}^{d}$,
$\{ \theta_{n}\}_{n=1}^{d}$,
$\{ \alpha_{n}\}_{n=1}^{d}$, satisfying the assumptions
of Proposition 3.6 and \st for $n$ odd $y_{n}$ is
a block of $(x_{k})_{k\in{\bf N}}$ while for $n$ even
$y_{n}$ is a block of $(w_{k})_{k\in{\bf N}}$.

\proof 
Let $j$ be given. We choose inductively a \sq
$\{ n_l\}_{l=1}^\infty\sbs\NN$ and vectors
$u_{l,A}, v_{l,A}\in X$, 
$u^*_{l,A}, v^*_{l,A}\in X^*$,  
$l=1,2,\ldots$, $A\sbs\{ 1,\ldots ,l-1\}$ \st
\begin{list}{(\alph{lista})} {\usecounter{lista}}
\item each $u_{l,A}$ is a block vector of 
$(x_k)_{k\in{\bf N}}$ and each $v_{l,A}$ is a block vector
of $(w_{k})_{k\in{\bf N}}$
\item For every $l=1,2,\ldots$ and every
$A\sbs\{ 1,\ldots ,l-1$ the vectors 
$u_{l,A}, u^*_{l,A}, v_{l,A},v^*_{l,A}$ are supported 
inside $(n_{l-1}, n_l]$.
\item Each $u_{l,A}$ is a $\left(\frac{1}{m^4_{2s}},2s\right)$-rapidly
increasing s.c.c., $u^*_{l,A}\in L_{2s}$
and $u^*_{l,A}\geq\frac{1}{4m_{2s}}$ where 
$2s>2j+1$ if $A=\emptyset$ and
$2s=\phi (u^*_{l_1,\emptyset},u^*_{l_2,A_1},\ldots ,
u^*_{l_k,A_{k-1}})$ if 
$A=\{ l_1<\cdots <l_k\}$ and 
$A_i=\{ l_1,\ldots ,l_i\}$,  $i=1,2,\ldots k-1 $.

The analogous relations hold for
$v_{l,A}, v^*_{l,A}$.
\end{list}

 The inductive construction
is straightforward.

Choose $F\sbs\{ n_l\}_{l=1}^\infty$,
$F=\{ n_{l_1},\ldots ,n_{l_k}\}$ \st a convex combination
$\sum_{n_l\in F}a_le_{n_l}$ is a 
$\left(\frac{1}{m^2_{2j+2}},2j+1\right)$-basic s.c.c.

For $i=1,\ldots ,k-1$, set
$ A_i=\{ l_1,\ldots ,l_i\}$.
Then it is easy to check that the \sq
\[ u_{l_1,\emptyset}, v_{l_2,A_1},\ldots ,v_{l_k,A_{k-1}}\]
and the corresponding one in $X^*$ have the desired properties.

\sprth
The following corollary is an  immediate consequence
of the previous proposition.

\sprt {\bf 3.13 Corollary.} The Banach space $X$ is Hereditarily
Indecomposable and hence it does not contain any unconditional basic
sequence.

\[{\bf REFERENCES }\]

\small
\begin{description}
\parsep =0pt
\itemsep =0pt
\parskip=0pt
\bb  {[Al-Ar]}]
  {\sc D. E. Alspach} and {\sc S. Argyros,} 
  {\sl Complexity of weakly null sequences,}
   Dissertationes Mathematicae, 321, 1992.
\bb {[Ar-D]}] {\sc S. Argyros} and {\sc I. Deliyanni,}
 {\sl Banach spaces of the type of Tsirelson,}
 preprint, 1992.
\bb {[B]}] {\sc S. F. Bellenot,}
{\sl Tsirelson superspaces and $\ell_p$,}
Journal of Funct. Analysis, 69, No 2, 1986, 207-228.
\bb{[C-S]}] {\sc P. G. Casazza} and {\sc T. Shura},
{\sl Tsirelson's space}, Lecture Notes in Math.
1363, Springer Verlag, 1989.
\bb{[F-J]}] {\sc T. Figiel} and {\sc W. B. Johnson},
{\sl A uniformly convex Banach space which contains no $\ell_p$},
Compositio Math. 29, 1974, 179-190.
\bb{[G1]}] {\sc W. T. Gowers,} 
{\sl A Banach space not containing $\ell_1$ or $c_0$ or
a reflexive subspace}, (preprint).
\bb{[G2]}] {\sc W. T. Gowers,} 
{\sl A hereditarily indecomposable space with an asymptotically
unconditional basis} (preprint).
\bb{[G3]}] {\sc W. T. Gowers,} 
{\sl A new dichotomy for Banach spaces} (preprint).
\bb {[G-M]}] {\sc W. T. Gowers} and {\sc B. Maurey,}
{\sl The unconditional basic \sq problem,}
Journal of AMS 6, 1993, 851-874.
\bb{[Ma]}] {\sc B. Maurey,}
{\sl A remark about distortion}, (preprint).
\bb{[Ma-Mi-To]}] {\sc M. Maurey, V. D. Milman} and
{\sc N. Tomczak-Jaegermann},
{\sl Asymptotic infinite dimensional theory of Banach spaces},
(preprint).
\bb{[Ma-R]}] {\sc B. Maurey} and {\sc H. Rosenthal},
{\sl Normalized weakly null sequences with no unconditional
subsequence},
Studia Math.  61, 1977, 77-98.
\bb{[Mi-To]}] {\sc V. D. Milman} and
{\sc N. Tomczak-Jaegermann},
{\sl Asymptotic $\ell_p$ spaces and bounded distortions},
Banach spaces, Contemp. Math. 144,
1993, 173-196.
\bb{[ O-Schl]}] {\sc E. Odell} and {\sc T. Schlumprecht},
The distortion problem, (preprint)
\bb{[ Schl]}] {\sc T. Schlumprecht}, {\sl An arbitrarily 
distortable Banach space}, Israel J. Math. 76,
1991, 81-95.
\bb{[Schr]}]{\sc J. Schreier,} {\sl Ein Gegenbeispiel zur
Theorie der sohwachen Konvergenz},
Studia Math. 2, 1930, 58-62. 
\bb {[T]}] {\sc B. S. Tsirelson,}
{\em Not every Banach space contains $\ell_{p}$ or
$c_{0}$,}
Funct. Anal. Appl. 8 (1974), p. 138-141.
\end{description}

\mbf e-mail: sargyros@atlas.uoa.ariadne-t.gr

\mbf deligia@talos.cc.uch.gr\ \ \ \ \ \ \ 
\end{document}